\documentclass[conference]{IEEEtran}
\IEEEoverridecommandlockouts
\usepackage[hidelinks]{hyperref}

\usepackage{multirow}
\usepackage[table,xcdraw]{xcolor}

\usepackage{cite}
\usepackage{caption}
\usepackage{subcaption}
\usepackage{algorithm, algorithmic}
\usepackage{amsmath,amssymb,graphicx,url} 
\newtheorem{assumption}{Assumption}
\newtheorem{theorem}{Theorem}
\newtheorem{remark}{Remark}
\newtheorem{lemma}{Lemma}

\DeclareMathOperator*{\argmin}{arg\,min}

\def\BibTeX{{\rm B\kern-.05em{\sc i\kern-.025em b}\kern-.08em
    T\kern-.1667em\lower.7ex\hbox{E}\kern-.125emX}}
\begin{document}

\title{Safe and Efficient Online Convex Optimization with Linear Budget Constraints and Partial Feedback
\thanks{
\IEEEauthorrefmark{1}Corresponding Author.}
}

\author{\IEEEauthorblockN{Shanqi Liu}
\IEEEauthorblockA{\textit{Research Institute of China Electronics Technology Group} \\
\textit{Shanghai, China}\\
liushanqi@stl.com.cn}
\and
\IEEEauthorblockN{Xin Liu\IEEEauthorrefmark{1}}
\IEEEauthorblockA{\textit{School of Information Science and Technology} \\
\textit{ShanghaiTech University, Shanghai, China}\\
liuxin7@shanghaitech.edu.cn}
}

\maketitle

\begin{abstract}
This paper studies online convex optimization with unknown linear budget constraints, where only the gradient information of the objective and the bandit feedback of constraint functions are observed. We propose a safe and efficient Lyapunov-optimization algorithm (SELO) that can achieve an $O(\sqrt{T})$ regret and zero cumulative constraint violation. The result also implies SELO achieves $O(\sqrt{T})$ regret when the budget is hard and not allowed to be violated.    
The proposed algorithm is computationally efficient as it resembles a primal-dual algorithm where the primal problem is an unconstrained, strongly convex and smooth problem, and the dual problem has a simple gradient-type update.  
The algorithm and theory are further justified in a simulated application of energy-efficient task processing in distributed data centers.  
\end{abstract}

\begin{IEEEkeywords}
safe online learning, bandit feedback, Lyapunov optimization, and optimistic/pessimistic design.
\end{IEEEkeywords}

\section{Introduction}
Online Convex Optimization (OCO) provides a versatile framework for studying online control and decision-making in dynamic and uncertain environments \cite{Sha_12,Haz_16, Fra_19}. Within this framework, a learner/controller continuously adapts its decisions to minimize a loss function or maximize a utility function while interacting with the environment in real-time. OCO has wide-ranging applications, including robotics  control \cite{WagCheJab_19, BarLamOli_21, SacBoo_22}, resource allocation in network systems \cite{CheLinGia_17, CheMokWan_17, AsgNee_20, YuNee_19, CaoBas_22, GuoCaoHe_23}, load balancing in server systems \cite{CaoZhaPor_18,HsuXuLin_22,JohKamKan_21}, online advertising \cite{McmHolScu_13, BalLuMir_20}, and personalized healthcare \cite{ZhaKeyEsm_20,TewSus_17}.

In OCO framework, the learner/controller chooses a decision $x_t$ in each round without prior knowledge of the convex loss function $f_t(\cdot)$. The objective is to minimize the total loss over $T$ rounds, given by $\min_{x\in \mathcal X}\sum_{t = 1}^T f_t(x_t),$ where the loss incurred by each decision $x_t$ is observed only after the decision has been executed. In real-world scenarios, these interactions are usually constrained by operational budgets. For instance, in a cloud computing system, the data center optimizes the throughput or delay performance within an energy budget; in online advertising, advertisers submit bids to maximize the click-through rates within a weekly or monthly budget. Similarly, inpatient onboarding, hospitals optimize patients' treatment by assigning them to proper medical units under medical resources constraints. 
In these applications, the operational constraints usually have a linear form, and the consumption matrix ($A_t$) is unknown, where  
only the partial feedback can be observed in practice (i.e., the gradient/values of loss functions $\nabla f(x_t)/f(x_t)$ and the values of resource consumption $A_t x_t$ of the decision point $x_t$).  
To model these applications, we study online convex optimization with linear budget constraints (OCOwLB) with partial feedback
\begin{align*}
    \min_{x_t \in \mathcal X}&  ~\sum_{t=1}^T f_t(x_t) \\
~~\hbox{s.t.}&~ \sum_{t=1}^T A_t x_t \leq B 
\end{align*}
where $\nabla f(x_t)/f(x_t)$ and $A_t x_t$ are revealed to the learner after the decision $x_t$ is taken. In OCOwLB, the learner aims to minimize the total loss with only partial feedback while satisfying the budget constraints or keeping the constraint violation minimal. 

One possible method to satisfy the budget constraints is to identify the underlying (expected) static constraints and project the decision $x_t$ into a safe (strictly feasible) set at each round \cite{ChaKal_22, ChaChaKal_23, AmaAliThr_19, AmaAliThr_20a, AmaAliThr_20b, AmaAliThr_22, HutCheAli_24}. However, such ``anytime safe projection'' methods may encounter three potential challenges when dealing with budget constraints: 1) they often require a substantial initial period to explore and learn the consumption matrix; 2) determining the ``correct'' safe constraint set based on an estimated consumption matrix is difficult and they are very likely to be overly conservative ensures safety but degrades performance; 3) the projection-based methods (e.g., projected online gradient descent) may require heavy computation because it is equivalent to solving a constrained quadratic optimization problem or could be almost extremely challenging (if not unsolvable) due to the non-convexity of the constraint set.
In this paper, we propose a safe and efficient Lyapunov-optimization (SELO) algorithm for OCOwLB with partial feedback that achieves strong theoretical results. Next, we summarize our main contributions. 

{\bf Main Contribution:} 
SELO is motivated by the primal-dual design and Lyapunov optimization framework. The primal decision modular incorporates the estimated loss function and the virtual queues that capture the over-used resource to balance the loss minimization and resource consumption. It is computationally efficient because the primal problem is (almost) unconstrained, strongly convex, and smooth.
The virtual queues $\{Q_t\}$ are designed to capture the cumulative overused budgets. In particular, we impose a ``pessimistic'' bonus term and ``pessimistic'' penalty term to guarantee the ``safety'' due to the initial exploration phase.    
SELO achieves $O(\sqrt{T})$ regret and zero constraint violation. 
We apply SELO to low-delay task processing in distributed data centers with energy budget constraints. Our experimental results demonstrate that SELO outperforms baseline algorithms and performs similarly against the optimal offline policy.

\section{Related work}
Safe online learning with partial/bandit feedback has been studied in the literature mainly from the perspective of constrained online convex optimization (COCO) \cite{MahJinYan_12, SunDeyKap_17, YuNeeWei_17, NeeYu_17, ChaSal_19, CaoZhaPoo_21, SadRauFaz_20, YuaLam_18, YiLiYan_21a, GuoLiuWei_22, LiuYanYin_23, AbhRah_24} and constrained bandit optimization (CBO) \cite{BadKleSli_13, BadKleSli_18, BadLanSli_14, AgrDev_16, ImmSanSch_22, SivZuoBan_22, BerCasCel_24b, CheCueLuo_20, Sin_24, SliSanFos_23, SliSanFos_24, GuoLiu_24, GuoLiu_25, PacGhaBar_21, CheGanSal_22, GanCheSal_24, GanGopSal_24}. 

{\bf COCO}: The work \cite{MahJinYan_12} initialized the studies of COCO and proposed an algorithm that achieves $O(\sqrt{T})$ regret and $O(T^{3/4})$ violation for static constraints and it is further extended to the setting where the constraints could be stochastic, time-varying or adversarial \cite{SunDeyKap_17, YuNeeWei_17, NeeYu_17, ChaSal_19, CaoZhaPoo_21, SadRauFaz_20}, where the representative work \cite{NeeYu_17} established  $O(\sqrt{T})$ regret and $O(\sqrt{T})$ constraint violation for stochastic and adversarial constraints when the Slater's condition holds. The work \cite{YuaLam_18} studied COCO under the strict metric of hard constraint violation where the violation cannot compensate and established $O(\sqrt{T})$ regret and $O(T^{3/4})$ hard violation for static constraints. The results is further improved to $O(\sqrt{T})$ regret and $O(T^{1/4})$ violation in \cite{YiLiYan_21a} and $O(\sqrt{T})$ regret and $O(1)$ violation in \cite{GuoLiuWei_22}. When the constraint is adversarial, the work \cite{GuoLiuWei_22} also proved $O(\sqrt{T})$ regret and $O(T^{3/4})$ hard violation and a recent remarkable work \cite{AbhRah_24} established $O(\sqrt{T})$ regret and $O(\sqrt{T}\log T)$ by designing a clever Lyapunov function integrated into adaptive online gradient descent. Unlike COCO with full information feedback, we only have partial feedback on the constraints.  

{\bf CBO}: The work \cite{BadKleSli_13, BadKleSli_18} initialized the studies of CBO by introducing the problem of Bandit with Knapsack and proposed the  algorithms (e.g., PrimalDualBwK) that achieve effective regret performance. BwK is further extended to various settings, including context (linear) bandit \cite{BadLanSli_14, AgrDev_16}, adversarial setting \cite{ImmSanSch_22, SivZuoBan_22, BerCasCel_24b}, fairness constraints \cite{CheCueLuo_20,Sin_24}, general constraints \cite{SliSanFos_23, GuoLiu_24}, and small budget setting \cite{GuoLiu_25}. These algorithms carefully balance the reward acquisition and resource consumption by using the idea of primal-dual optimization more or less. Note another line of CBO is to satisfy the constraints for every round \cite{PacGhaBar_21, CheGanSal_22, GanCheSal_24, GanGopSal_24}, where these papers usually construct the safe or feasibility region and then choose the decision by solving a (linear) optimization problem for every round. Unlike these literature, we focus on the continuous (convex) decision space.

\section{Model}\label{sec: model}
In this section, we formally define OCOwLB. Consider the following online learning problem: At each round $t \in [T]$, the learner makes decision $\mathbf x_t \in \mathcal X$ and then observes the gradient or value of loss function $\nabla f_t(\mathbf x_t)/f_t(\mathbf x_t)$ and budget consumption $\mathbf o_t = \mathbf A_t \mathbf x_t$ after the decision $x_t$ is executed.  
The feedback differs from the classical OCO where only the bandit feedback of the constraint function can be observed \cite{Sha_12, Haz_16, Fra_19}. We assume $\{f_t(\cdot)\}_t$ is convex and  arbitrarily time-varying and $\{\mathbf A_t\}_t$ is an i.i.d. random sequence with $\mathbb E[\mathbf A_t] = \mathbf A, \forall t \in [T]$. 
We study two types of (soft and hard) budget constraints with initial budget $\mathbf b_T$: 1) in soft budget constraints, the slightly overused budgets are allowed and the interaction can continue even when the budget is used up; 2) in hard budget constraints, the violation is not allowed, and the interaction terminates when the budget is used up.

The goal of the learner is to generate a decision sequence $\{\mathbf x_1,\mathbf x_2,\ldots,\mathbf x_{T-1},\mathbf x_T\}$ to minimize the total loss $\sum_{t=1}^T f_t(\mathbf x_t)$ by best using the resource budget. To quantify the performance of an OCOwLB, we compare it with a static baseline, called the best-fixed decision in hindsight, which is the solution to the following offline OCOwLB such that
\begin{align}
    \min_{\mathbf x \in \mathcal X}&  ~\sum_{t=1}^T f_t(\mathbf x)
    \label{offline coco intro}\\
~\hbox{s.t.}&~ \mathbf A \mathbf x \leq \mathbf b. \label{offline coco cons} 
\end{align}
where $\mathcal X$ is a ``simple'' convex set (e.g., positive quadrant or probability simplex) and $\mathbf b = \mathbf b_T/T.$ 
The optimal solution $\mathbf x^*$ to offline OCOwLB \eqref{offline coco intro}--\eqref{offline coco cons} is called the best single
decision in hindsight, a widely used baseline for the online learning literature \cite{Sha_12, Haz_16, Fra_19}. 
Based on offline OCOwLB, the regret and constraint violation (for soft budget) of an online algorithm are defined below
\begin{align*}
    \mathcal R(T) =&~ \mathbb E[\sum_{t=1}^T (f_t(\mathbf x_t) - f_t(\mathbf x^*))], \\
    \mathcal V(T) =&~ \max(\mathbb E[\sum_{t=1}^T \mathbf A_t \mathbf x_t - \mathbf b_T], 0). \nonumber
\end{align*}
The expectation is taken over the randomness of the algorithm and environment. 
Note that for the hard budget constraints, 
we only measure the regret of an online learning algorithm because no constraint violation exists. In this paper, we show they are equivalent under our algorithm and analysis because the regret and constraint violation can trade off each other (once we have the anytime high probability guarantee on the constraint violation). The details can be found in Section \ref{sec: theory}.   

\section{Safe and efficient Lyapunov \\ optimization algorithm} \label{sec: alg}
We present a safe and efficient Lyapunov optimization algorithm (SELO) for OCOwLB with bandit feedback. In SELO, we first need to estimate the loss functions (via the gradient feedback) and the safe region via bandit feedback. For the loss function, 
we define the estimated loss function
\begin{align}
\hat f_{t}(\mathbf x) =&~ f_{t-1}(\mathbf x_{t-1})+ \langle \nabla f_{t-1}(\mathbf x_{t-1}),  \mathbf x- \mathbf x_{t-1} \rangle. \label{eq: est f} 
\end{align}
For the unknown consumption matrix, we leverage the optimistic/pessimistic principle in \cite{AueCesFis_02, AbbPalSze_11, LatSze_20} to estimate $\mathbf A$ from the bandit feedback. At the end of round $t-1,$ we estimate the $\mathbf A$ with the historical observation $\{\mathbf x_s, \mathbf o_s\}_{s=1}^{t-1}$ using the regularized least square regression 
\begin{align}
\mathbf{\bar A}_t = \mathbf \Sigma_t^{-1} \sum_{s=1}^{t-1} \mathbf o_s \mathbf x_{s},~~\mathbf \Sigma_t = \mathbf I + \sum_{s=1}^{t-1} \mathbf x_{s}\mathbf x^{\dag}_{s}. \label{eq: est A} 
\end{align}
Based on the optimistic principle in \cite{AueCesFis_02, AbbPalSze_11, LatSze_20}, we can construct the pessimistic consumption for any $\mathbf x \in \mathcal X$ such that $$\mathbf{\bar A}_t \mathbf x + \alpha\|\mathbf x\|_{\mathbf \Sigma_t^{-1}} \geq \mathbf A \mathbf x,$$ where $\alpha=\Theta(\sqrt{\log T} + \sqrt{\boldsymbol{\lambda}})$ is the radius and the inequality holds with a high probability \cite{AueCesFis_02, AbbPalSze_11, LatSze_20}. 
Given the pessimistic estimation of resource consumption, we  define the pessimistic constraint functions 
\begin{align}
\hat g_t(\mathbf x) =&~ \mathbf{\bar A}_t \mathbf x + \alpha \|\mathbf x\|_{\mathbf \Sigma_t^{-1}} - \mathbf b.  \label{eq: est g}
\end{align}
Now we are ready to formally present our proposed SELO algorithm based on the introduced estimators.
\begin{algorithm}[H]
\caption{SELO Algorithm}
\begin{algorithmic}[1]
\STATE {\bf Initialization:} $Q_{0} = 0, T_0, \eta, \xi, \alpha$ and $V.$ 
\FOR{$t=1,\cdots, T_0$}
\STATE {\bf (Short) pure exploration phase:} generate a (scaled) Gaussian
random vector $\mathbf x_t \in \mathcal X.$
\ENDFOR
\FOR{$t=T_0+1,\cdots, T$}
    \STATE {\bf Safe and efficient decision:}  
    \begin{align*}
    \mathbf \mathbf x_{t} = \argmin_{\mathbf x \in  \mathcal X} ~  V \hat f_t(\mathbf x) + \langle \mathbf Q_t, \hat g_t(\mathbf x) \rangle +   \frac{1}{2\eta} \|\mathbf x-\mathbf x_{t-1}\|^2. 
    \end{align*}
    \STATE {\bf Double pessimistic budget pacing:} 
    \begin{align*}
     &\mathbf Q_{t+1} = \max(\mathbf Q_{t} + \hat g_t(\mathbf x_{t}) + \xi \mathbf 1, 0). 
     \end{align*}
    \STATE {\bf Observe the feedback $\{f_{t}(\mathbf x_t), \nabla f_{t}(\mathbf x_t), \mathbf o_t\}$ and update the estimation according to \eqref{eq: est f}--\eqref{eq: est g}}. 
\ENDFOR
\end{algorithmic} \label{alg}
\end{algorithm}
We explain the intuition behind SELO from the perspective of primal-dual optimization in the offline problem \eqref{offline coco intro}--\eqref{offline coco cons}. The Lagrange function of the offline optimization problem in \eqref{offline coco intro}--\eqref{offline coco cons} can be written as   $$L(\boldsymbol{\lambda}, \mathbf x) := \sum_{t=1}^T L_t(\boldsymbol{\lambda}, \mathbf x) := \sum_{t=1}^T (f_t(\mathbf x) + \langle \boldsymbol{\lambda}, \mathbf A \mathbf x \rangle),$$
where the dual variables $\boldsymbol{\boldsymbol{\lambda}}$ are associated with the constraints in \eqref{offline coco cons}. 
Since we have no prior knowledge of $f_{t}(\cdot)$ when making a decision $\mathbf x_{t},$ we replace it with its estimator in \eqref{eq: est f}. 
Similarly, we utilize the pessimistic $\hat g_t(\mathbf x)$ to replace $\mathbf A \mathbf x - \mathbf b.$ 
For the key dual variable, we design a virtual queue to mimic the optimal dual $\boldsymbol{\lambda}^*$, which incorporates the pessimistic design of $\hat g_t(\mathbf x_{t})$ and the pessimistic factor of $\xi$. Specifically, the virtual queues $\mathbf Q_t$ are to monitor the cumulative consumption of the budget. An increase in these virtual queues indicates potential overuse of resources, prompting more conservative decision-making. The virtual queues can be viewed as scaled dual variables within the primal--dual framework described in \cite{Nee_10, SriYin_14}.
Furthermore, the primal decision can be interpreted as minimizing ``loss + Lyapunov drift'' as $\frac{\mathbf Q_{t+1}^2-\mathbf Q_{t}^2}{2} \approx \langle \mathbf Q_t, \hat g_t(\mathbf x) \rangle$''.
Therefore, we can use Lyapunov drift analysis to establish the upper bound of the virtual queues, which can then be translated into the violation bound.

Finally, we comment that our algorithm is computationally efficient as it only needs to solve an ``almost'' unconstrained optimization problem (if $\mathcal X$ is a simple set like the boundary constraints). When $\mathcal X \in \mathbb R^d$ is general convex, the optimization problem in line-$6$ can be formulated as a quadratically constrained quadratic program (QCQP). Therefore, by using interior point methods, the time complexity is $O(d^{3.5} \log(1/\epsilon))$ to establish an $\epsilon-$accuracy solution and the space complexity is $O(d^2).$
We can show the surrogate function of $V\hat f_t(\mathbf x) + \langle \mathbf Q_t, \hat g_t(\mathbf x) \rangle +   \frac{1}{2\eta} \|\mathbf x-\mathbf x_{t-1}\|^2$ is both strongly convex and smooth by showing its hessian matrix is positive definite \cite{BoyVan_04}.
Therefore, we could use the classical gradient descent algorithm to find the minimizer with a linear and dimension-free converge rate \cite{SchRouBac_11}. It is worth to be mentioned that our algorithm is much more efficient than the projection online gradient descent algorithm \cite{Zin_03, MahJinYan_12}, and the online Frank-Wolfe algorithm \cite{HazKal_12}. 

\section{Theoretical Results}\label{sec: theory}
We next analyze the regret  $\mathcal R(T)$ and violation $\mathcal V(T)$ of SELO based on the following standard technical assumptions on the feasible set, loss, and constraint functions.  

\begin{assumption}
The feasible set $\mathcal X$ is convex with diameter $D$ such that $\|\mathbf x - \mathbf x'\| \leq D, \forall \mathbf x, \mathbf x' \in \mathcal X$. \label{assumption: set}
\end{assumption}


\begin{assumption}
The loss function is convex and Lipschitz continuous with Lipschitz constant $F$ such that $|f_t(\mathbf x) - f_t(\mathbf x')| \leq F \|\mathbf x - \mathbf x'\|, \forall \mathbf x, \mathbf x' \in \mathcal X, \forall t.$ The consumption matrix $\{\mathbf A_t\}_1^{T}$ is i.i.d. generated with $\mathbb E[\mathbf A_t] = \mathbf A, \forall t,$ where each element is bounded $\mathbf A_{t}(i,j) \in [0,1], \forall i, j, t$. The average budget $\mathbf b \in [0, 1].$ \label{assumption: fun}
\end{assumption}

\begin{assumption}
There exists a positive value $\beta$ and strictly feasible $\mathbf x \in \mathcal X$ and such that $\mathbf A \mathbf x - \mathbf b  \leq -\beta \mathbf 1.$ 
\label{assumption: slater}
\end{assumption}

Note Assumptions \ref{assumption: set} and \ref{assumption: fun} imply the feasible set $\mathcal X$ and the functions are bounded. Assumption \ref{assumption: slater} is Slater's condition and common in the (online) optimization literature \cite{BoyVan_04, SchRouBac_11}. 
With these assumptions, we are ready to present our main theoretical results of SELO w.r.t. the budget constraints are soft or hard. 
\begin{theorem}\label{thm: main}
Lett  $V = \sqrt{T},$ $\eta = 1/T,$ $\xi = \log^2 T/\sqrt{T}$ and $T_0=\log T/\beta.$ Under Assumptions \ref{assumption: set}--\ref{assumption: slater}, SELO algorithm achieves the following theoretical results for OCOwLB

\noindent 1) under soft budget constraints such that
\begin{align*}
\mathcal R(T) = \Tilde{O}(\sqrt{T}),
~\mathcal V(T) = O(1). 
\end{align*} 

\noindent 2) under hard budget constraints such that
\begin{align*}
\mathcal R(T) = \Tilde{O}(\sqrt{T}).
\end{align*} 
\end{theorem}

\begin{remark} \label{mark}
Theorem \ref{thm: main} establishes the  $\Tilde{O}(\sqrt{T})$ regret and $O(1)$ violation for the soft budget constraint and $\Tilde{O}(\sqrt{T})$ regret for the hard budget constraint. These results recover $\Tilde{O}(\sqrt{T})$ regret for the classical OCO without budget constraints. This is optimal up to a $\log$ factor because the lower bound of regret in OCO is $\Omega(\sqrt{T})$ \cite{Sha_12,Haz_16, Fra_19}. The results are better than ``anytime safe projection'' methods because the work \cite{ChaKal_22} incurs $\Tilde{O}(T^{\frac{2}{3}})$.
When the only bandit feedback of loss functions is known (i.e., only the value $f_t(\mathbf x_t)$ is observed), one could leverage either a two-point or one-point gradient estimation method into
SELO algorithm \cite[Chapter 6]{Haz_16}. 
\end{remark}

To establish these strong results, we need to carefully analyze the reward acquisition and budget consumption processes. The key is 1) identifying an optimal but slightly conservative baseline/set and 2) the upper bound of the virtual queues. 
For 1), we establish a safe and close-to-optimal set where the safety margin is diminishing as the estimation errors of the consumption matrix; The accumulated degraded regret performance can be bounded by $\Tilde{O}(\sqrt{T});$ 
For 2), our analysis regards the virtual queue update as a Markovian process and leverages the multiple-step Lyapunov-drift analysis to establish the high probability bound on the virtual queues (i.e., the over-consumed budgets). The bounded virtual queue will be translated into the soft constraint violation and the additional regret due to inappropriate early stopping under the hard budget constraint. 

As discussed in Section \ref{sec: alg},  our algorithm can be interpreted as minimizing ``regret + Lyapunov drift''. 
Let the Lyapunov function and its drift be 
\begin{align*}
    L_t = ~\frac{1}{2}\|\mathbf {Q}_{t}\|^2_2, ~~\Delta_t = ~L_{t+1} - L_{t}.
\end{align*}
Let $\mathbb D(\mathbf x,\mathbf y,\mathbf z) = \|\mathbf x-\mathbf y\|^2 - \|\mathbf x-\mathbf z\|^2.$
Define $\mathbb E_{\mathcal H_t}[\cdot] = \mathbb E[\cdot | \mathcal H_t],$ where $\mathcal H_t =\{\mathbf x_s, \nabla f_s, \mathbf o_s, \mathbf {Q_s}\}_{s=1}^{t-1}.$ 

Further, we introduce the important concept $\epsilon_t-$tight feasible set to the offline problem \eqref{offline coco intro}-\eqref{offline coco cons} defined by $$\mathcal X_{\epsilon_t} = \{\mathbf x \in \mathcal X ~|~ \mathbf A \mathbf x \leq \mathbf b - \epsilon_t \mathbf 1\},$$ where $\epsilon_t = 2(\xi + 2 \alpha\|\mathbf x\|_{\mathbf \Sigma_t^{-1}}).$ We establish the following key lemma that bridges the one-step regret and Lyapunov drift. 
\begin{lemma}\label{lem: key} 
Under Assumptions \ref{assumption: set}--\ref{assumption: slater}, SELO in Algorithm \ref{alg} establishes that for any  $\epsilon_t-$tight feasible solution $\mathbf x \in \mathcal X_{\epsilon_t}$ that
\begin{align}
&~~~~\mathbb E_{\mathcal H_t} \left[V(f_t(\mathbf x_t) - f_t(\mathbf x)) + \Delta_t \right] \label{eq: key}\\
&\leq \mathbb E_{\mathcal H_t} \left[\langle \mathbf Q_t, \mathbf A \mathbf x - \mathbf b\rangle + \frac{\mathbb D(\mathbf x, \mathbf x_{t},\mathbf x_{t+1})}{2\eta} + F^2\right]. \nonumber
\end{align} 
\end{lemma}
The key lemma establishes the bound of ``one-step regret + Lyapunov drift'' in \eqref{eq: key} and bridges the analysis to bound both regret and the virtual queue (i.e., constraint violation).
The upper bound includes the key cross-term, the proximal bias $\mathbb D(\mathbf x, \mathbf x_t, \mathbf x_{t+1})$, trade-off factor $\xi.$
We choose the parameters $\{V, \eta, \xi\}$ to minimize the upper bound.
Next, the key is to quantify the cross term in \eqref{eq: key}, which is related to the way we choose the safe and optimal benchmark $\mathbf x \in \mathcal X_{\epsilon_t}$. 
Therefore, we need to carefully choose a ``safe and optimal'' baseline to analyze the regret and constraint violation. 
Next, we sketch the proof of Theorem \ref{thm: main} based on the key Lemma \ref{lem: key}. 

{\noindent \bf Regret analysis:} let $\mathbf x = \mathbf x_{\epsilon_t}^*$ be the offline optimal solution w.r.t. $\mathbf x \in \mathcal X_{\epsilon_t}$ and we immediately have $\mathbf A \mathbf  x_{\epsilon_t}^* - b \mathbf 1 \leq 0.$ Therefore, the cross term becomes negative in \eqref{eq: key}, and we take the telescope summation from $t=1$ to $T$ such that
\begin{align}
\sum_{t=1}^T(f_{t}(\mathbf x_{t}) - f_t(\mathbf x^*_{\epsilon_t}))
= O\left(\frac{D^2}{V\eta} + T \xi + \frac{T}{V}\right). \label{eq: key regret}
\end{align}
Recall $V = \sqrt{T},$ $\eta = T,$ $\xi = \log^2 T/\sqrt{T}$ and we have \eqref{eq: key regret} bounded by $\Tilde{O}(\sqrt{T}).$ 
The original regret $\mathcal R(T)$ is also $\Tilde{O}(\sqrt{T})$ because the cumulative error $\sum_{t=1}^T \epsilon_t = \Tilde{O}(\sqrt{T})$ according to \cite[Theorem 2]{AbbPalSze_11}. 

{\noindent \bf Violation analysis:} we choose a strictly feasible $\mathbf x$ such that $\mathbf A \mathbf x - \mathbf b \leq -\beta\mathbf 1.$ It implies that 
\begin{align}
&\mathbb E_{\mathcal H_t} \left[\Delta_t \right] \label{eq: key Q}\\
\leq& -\beta\|\mathbf Q_t\| + O\left(\frac{\mathbb D(\mathbf x,\mathbf x_{t}, \mathbf x_{t+1})}{2\eta} + V (FD + \xi) \right). \nonumber
\end{align} 
Intuitively, if the virtual queue process $\{\mathbf Q_t\}$ already reaches
the steady state, i.e., the drift is zero $\mathbb E[\Delta(t)] = 0$ and $\mathbb E[D(\mathbf x,\mathbf x_{t}, \mathbf x_{t+1})] = 0.$ We immediately establish its upper bound of $\mathbb E[Q_t] = O(\frac{V(FD+\xi)}{\beta}) = O(\sqrt{T}).$ This is an ideal constraint violation bound.   
Motivated by this intuition, we analyze its finite time performance by constructing a multi-step Lyapunov drift analysis, which is formally justified by the Forster-Lyapunov lemma \cite{YuNeeWei_17, LiuLiShiYin_21} and achieve anytime high probability bound 
\begin{align}
\mathbb P(\|\mathbf Q_t\| \leq O(\sqrt{T}\log T)) \leq 1 - T^{-2}, ~\forall t \in [T].    \label{eq: Q bound}
\end{align}
The bounded virtual queue immediately implies a zero violation because for a large $T$ such that $T \xi = \Omega(\sqrt{T}\log^2 T).$ 
$$\mathcal V(T) \leq (\mathbb E[\|\mathbf Q_T\|_1] - T \xi)^{+}.$$

\noindent {\bf Regret analysis under hard budget constraints:}
Finally, for the setting with hard budget constraints, $\|\mathbf Q_t\|$ is bounded anytime by $O(\sqrt{T}\log T)$ in \eqref{eq: Q bound}, which means the interaction could stop early due to the over-consumed budget resource. This amount of violation can be translated into the corresponding regret according to the budget pacing such that
\begin{align*}
 \sum_{t=1}^T  g_t(\mathbf x_{t}) \leq \mathbf Q_{T+1} + \sum_{t=1}^T ( g_t(\mathbf x_{t}) -  \hat g_t(\mathbf x_{t})). 
\end{align*}
Therefore, we have at most $\Tilde{O}(\sqrt{T})$ additional regret due to the possible early stopping due to the hard budget constraint. 

As discussed in theoretical results, the key is to prove the key Lemma ``regret + Lyapunov drift". We will focus on proving this key lemma in the next section. 


\section{Theoretical Analysis of the key ``Regret + Lyapunov Drift'' lemma}\label{sec: theory}
To prove the key lemma, we first introduce an important Lemma from \cite{GonMar_93, YuNee_20, Lin_22}, which provides a useful upper bound on the optimal value of the strongly convex function. The proof is straightforward by definitions of the strongly convex function and the first-order optimality condition of the convex function. 
\begin{lemma}
Let $\mathcal X$ be a convex set. Let $h: \mathcal X \to \mathcal R$ be $\alpha$-strongly convex function on $\mathcal X$ and $x_{opt}$ be an optimal solution of $h,$ i.e., $x_{opt} = \argmin_{x\in\mathcal X} h(x).$ Then, $h(x_{opt}) \leq h(x) - \frac{\alpha}{2}\|x - x_{opt}\|^2$ for any $x\in \mathcal X.$ 
\label{lemma:tool}
\end{lemma}

We then provide the detailed proof for Lemma \ref{lem: key}. 
Recall the decision in Algorithm \ref{alg}, we have
\begin{align*}
    \mathbf \mathbf x_{t} = \argmin_{\mathbf x \in  \mathcal X} ~  V\hat f_t(\mathbf x) + \langle \mathbf Q_t, \hat g_t(\mathbf x) \rangle +   \frac{1}{2\eta} \|\mathbf x-\mathbf x_{t-1}\|^2. 
\end{align*}
According to Lemma \ref{lemma:tool}, we have
\begin{align*}
&V\hat f_t(\mathbf x_t) + \langle \mathbf Q_t, \hat g_t(\mathbf x_t) \rangle +   \frac{1}{2\eta} \|\mathbf x_t-\mathbf x_{t-1}\|^2 \\
\leq& V\hat f_t(\mathbf x) + \langle \mathbf Q_t, \hat g_t(\mathbf x) \rangle +   \frac{1}{2\eta} \|\mathbf x-\mathbf x_{t-1}\|^2 - \frac{1}{2\eta} \|\mathbf x-\mathbf x_{t}\|^2
\end{align*}
Recall the definition of $\hat f$ in \eqref{eq: est f} $$\hat f_{t}(\mathbf x) = f_{t-1}(\mathbf x_{t-1})+ \langle \nabla f_{t-1}(\mathbf x_{t-1}),  \mathbf x-\mathbf x_{t-1} \rangle,$$ we have for any $\mathbf x \in \mathcal X$ such that $$\hat f_{t}(\mathbf x) \leq f_{t-1}(\mathbf x)$$ due to the convexity of function $f(\mathbf x).$ 
It implies that  
\begin{align*}
&f_{t-1}(\mathbf x_{t-1}) - f_{t-1}(\mathbf x) - \frac{\eta}{2} \|\nabla f_{t-1}(\mathbf x_{t-1})\|^2\\
\leq& \hat f_t(\mathbf x_t) - \hat f_t(\mathbf x) +   \frac{1}{2\eta} \|\mathbf x_t-\mathbf x_{t-1}\|^2. 
\end{align*}
Therefore, recall the definition of $D(\cdot, \cdot, \cdot)$ and the bounded assumption in Assumption \ref{assumption: fun}, we have  
\begin{align*}
&V(f_{t-1}(\mathbf x_{t-1}) - f_{t-1}(\mathbf x)) + \langle \mathbf Q_t, \hat g_t(\mathbf x_t) \rangle \\
\leq& \langle \mathbf Q_t, \hat g_t(\mathbf x) \rangle +   \frac{D(\mathbf x, \mathbf x_{t},\mathbf x_{t+1})}{2\eta} +   \frac{\eta V^2F^2}{2}. 
\end{align*}
Combine with the virtual queue update,
\begin{align*}
 &\mathbf Q_{t+1} = \max(\mathbf Q_{t} + \hat g_t(\mathbf x_{t}) + \xi \mathbf 1, 0),
 \end{align*}
we have the one-step drift
\begin{align*}
\frac{1}{2} \|\mathbf Q_{t+1}\|^2 - \frac{1}{2}\|\mathbf Q_t\|^2 \leq \langle \mathbf Q_t, \hat g_t(\mathbf x_t ) + \xi \mathbf 1 \rangle,
\end{align*}
because $(\max(x,0))^2\leq x^2$ holds. 
It further implies that 
\begin{align*}
&V(f_{t-1}(\mathbf x_{t-1}) - f_{t-1}(\mathbf x)) + \frac{1}{2} \|\mathbf Q_{t+1}\|^2 - \frac{1}{2}\|\mathbf Q_t\|^2  \\
\leq& \langle \mathbf Q_t, \hat g_t(\mathbf x) + \xi \mathbf 1 \rangle +   \frac{\mathbb D(\mathbf x, \mathbf x_{t},\mathbf x_{t+1})}{2\eta} +   \frac{\eta V^2F^2}{2}. 
\end{align*}
Now, let's study the key cross term 
\begin{align*}
&\langle \mathbf Q_t, \hat g_t(\mathbf x) + \xi \mathbf 1 \rangle  \\
=& \langle \mathbf Q_t, \mathbf{\bar A}_t \mathbf x + (\alpha \|\mathbf x\|_{\mathbf \Sigma_t^{-1}}  + \xi) \mathbf 1  - \mathbf b  \rangle \\
=& \langle \mathbf Q_t, \mathbf{A} \mathbf x - \mathbf b  + (\mathbf{\bar A}_t - \mathbf{A}) \mathbf x +  (\alpha \|\mathbf x\|_{\mathbf \Sigma_t^{-1}}  + \xi) \mathbf 1   \rangle \\
\leq& \langle \mathbf Q_t, \mathbf{A} \mathbf x - \mathbf b + 2(\alpha \|\mathbf x\|_{\mathbf \Sigma_t^{-1}}  + \xi) \mathbf 1   \rangle 
\end{align*}
where the last inequality holds according to the high probability event of the upper confidence bound in \cite{AbbPalSze_11} such that
$$\mathbb P((\mathbf{\bar A}_t - \mathbf{A}) \mathbf x \geq (\alpha \|\mathbf x\|_{\mathbf \Sigma_t^{-1}}  + \xi) \mathbf 1  ) \leq 1-T^{-3}.$$
Recall the set $\mathcal X_{\epsilon_t}$ where $\epsilon_t = 2(\xi + 2 \alpha\|\mathbf x\|_{\mathbf \Sigma_t^{-1}}),$ we have the key inequality for any $\mathbf x \in \mathcal X_{\epsilon_t}$ such that 
\begin{align*}
\langle \mathbf Q_t, \hat g_t(\mathbf x) + \xi \mathbf 1 \rangle  
\leq \langle \mathbf Q_t, \mathbf{A} \mathbf x - \mathbf b  \rangle 
\end{align*}
Finally, we have 
\begin{align}
&~~~~\mathbb E_{\mathcal H_t} \left[V(f(\mathbf x_t) - f(\mathbf x)) + \frac{1}{2} \|\mathbf Q_{t+1}\|^2 - \frac{1}{2}\|\mathbf Q_t\|^2 \right] \label{eq: key}\nonumber\\
&\leq \mathbb E_{\mathcal H_t} \left[\langle \mathbf Q_t, \mathbf A \mathbf x - \mathbf b\rangle + \frac{\mathbb D(\mathbf x, \mathbf x_{t},\mathbf x_{t+1})}{2\eta} +  F^2\right], \nonumber
\end{align} 
Recall the definition of $\Delta_t = \frac{1}{2} \|\mathbf Q_{t+1}\|^2 - \frac{1}{2}\|\mathbf Q_t\|^2$ and the key Lemma \ref{lem: key} is proved. 

\section{Experiments}
\label{sec: experiments}
In this section, we evaluate the performance of SELO in energy resource allocation to optimize the delay performance in distributed data centers within the energy budget constraints. 
We considered a distributed data center with server clusters located in different regions where the incoming jobs arrive at each region and for the service. The service capability of a cluster is a function w.r.t. its energy consumption and the energy prices vary across locations and times. The goal is to minimize the queueing delay of jobs while guaranteeing the energy cost within the budget. This problem can be formulated as an OCOwLB and solved by our algorithm.

Specifically, we considered regions with $r=10$, and the duration of each round was one hour. At round $t,$ the average number of arrival tasks is $\boldsymbol{\lambda}_t \in \mathbb R^{10}$; the service rate is $\boldsymbol{\mu}_t \in \mathbb R^{10}$ and $\mathbf x_t\in \mathbb R^{10}$ be the fraction of service capacity on at the round $t.$ We model the data centers in each region as M/M/1 queueing system \cite{Har_13}, and the delay is approximated by $d_t(\mathbf x_{t,i}) = 1/(\mathbf e_i + \mathbf x_{t,i} \boldsymbol{\mu}_{t,i} - \boldsymbol{\lambda}_{t,i})$ with $\mathbf e_i$ being the base service capacity at region $i$. The loss function (system-wise delay) at round $t$ is $f_t(\mathbf x_t) = \sum_{i=1}d_i (\mathbf x_{t,i}).$  
Let $\mathbf A_t \in \mathbb R^{10}$ are the energy price vector at round $t$ and $\langle \mathbf A_t, \mathbf  x_t \rangle$ is the energy costs at round $t.$ 
In the experiment, we extract and calibrate the average electricity price from the trace 
in New York City from NYISO \cite{NewYorISO}. 
We set the unit base service at each region $\mathbf e_i=1, \forall i.$ The traffic $\{\boldsymbol{\lambda}_t\}_{t=1}^T$ is generated according to Gaussian random variables where the (normalized) mean values are time-varying to mimic realistic daily periodic patterns. The service $\{\boldsymbol{\mu}_t\}_{t=1}^T$ is Gaussian random variables draw from $\mathcal N(5, 0.5).$ The average budget $b = 0.75.$
\begin{figure} 
\centering
\begin{subfigure}{0.43\textwidth}
    \includegraphics[width=\textwidth]{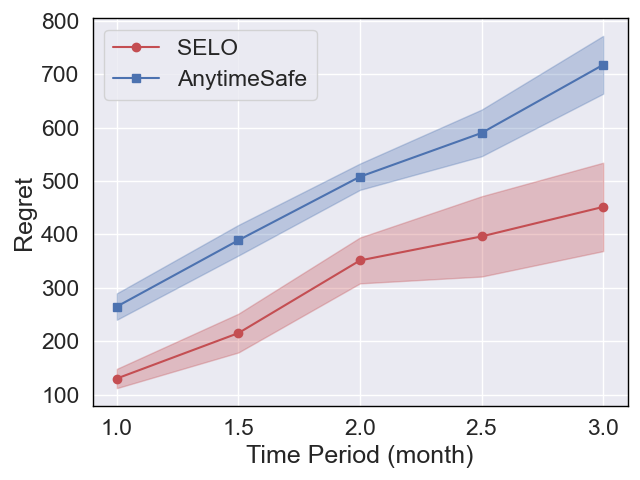}
\end{subfigure}
\hfill
\begin{subfigure}{0.43\textwidth}
    \includegraphics[width=\textwidth]{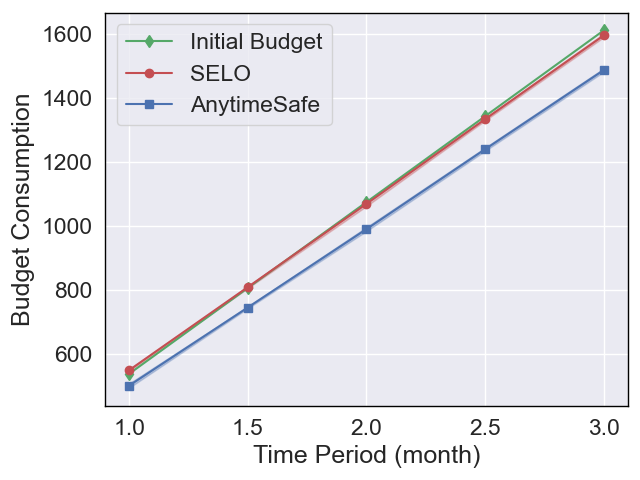}
\end{subfigure}
\caption{Performance comparison: SELO v.s. ``AnytimeSafe".}
\label{fig: alg}
\end{figure}

We compared SELO with the ``AnytimeSafe" algorithm from \cite{ChaKal_22}, best-tuning its initial period while keeping other hyperparameters the same as ours. We solve the offline problem in \eqref{offline coco intro}--\eqref{offline coco cons}, which serves as the optimal baseline for comparing and analyzing regret. We tested these algorithms over a range of one to three months. 
The results shown in Figure \ref{fig: alg} demonstrate that 
both SELO and ``AnytimeSafe" algorithms can achieve safe and zero constraint violation performance. However, 
SELO achieves superior delay performance by effectively and safely utilizing the resource budget, whereas ``AnytimeSafe" tends to be over-conservative to guarantee safety and might incur a large regret.


\section{Conclusions}
In this paper, we studied online convex optimization with unknown linear constraints. We proposed a safe and efficient Lyapunov optimization (resembles a primal-dual algorithm) with a dedicated design of double pessimistic budget pacing, and we proved it achieves strong regret and constraint violation performance by utilizing a multiple-step Lyapunov drift analysis. 
The simulated experiments on energy-efficient task load balancing in distributed data centers justified our algorithm and theoretical analysis.


\bibliographystyle{ieeetr}
\bibliography{ref.bib}

\begin{thebibliography}{10}

\bibitem{Sha_12}
S.~Shalev-Shwartz, ``Online learning and online convex optimization,'' {\em Found. Trends Mach. Learn.}, vol.~4, no.~2, p.~107–194, 2012.

\bibitem{Haz_16}
E.~Hazan, ``Introduction to online convex optimization,'' {\em Foundations and Trends{\textregistered} in Optimization}, vol.~2, no.~3-4, pp.~157--325, 2016.

\bibitem{Fra_19}
F.~Orabona, ``A modern introduction to online learning,'' {\em arXiv preprint arXiv:1912.13213}, 2019.

\bibitem{WagCheJab_19}
N.~Wagener, C.-A. Cheng, J.~Sacks, and B.~Boots, ``An online learning approach to model predictive control,'' 2019.

\bibitem{BarLamOli_21}
{Barcelos, Lucas and Lambert, Alexander and Oliveira, Rafael and Borges, Paulo and Boots, Byron and Ramos, Fabio}, ``{Dual Online Stein Variational Inference for Control and Dynamics},'' in {\em Robotics: Science and Systems ({R:SS})}, 2021.

\bibitem{SacBoo_22}
{Sacks, Jacob and Boots, Byron}, ``{Learning to Optimize in Model Predictive Control},'' in {\em {IEEE} International Conference on Robotics and Automation ({ICRA})}, 2022.

\bibitem{CheLinGia_17}
T.~Chen, Q.~Ling, and G.~B. Giannakis, ``An online convex optimization approach to proactive network resource allocation,'' {\em IEEE Transactions on Signal Processing}, vol.~65, no.~24, pp.~6350--6364, 2017.

\bibitem{CheMokWan_17}
T.~Chen, A.~Mokhtari, X.~Wang, A.~Ribeiro, and G.~B. Giannakis, ``Stochastic averaging for constrained optimization with application to online resource allocation,'' {\em IEEE Transactions on Signal Processing}, 2017.

\bibitem{AsgNee_20}
K.~Asgari and M.~J. Neely, ``Bregman-style online convex optimization with energyharvesting constraints,'' {\em Proc. ACM Meas. Anal. Comput. Syst. (POMACS)}, 2020.

\bibitem{YuNee_19}
H.~Yu and M.~J. Neely, ``Learning-aided optimization for energy-harvesting devices with outdated state information,'' {\em IEEE/ACM Transactions on Networking}, vol.~27, no.~4, pp.~1501--1514, 2019.

\bibitem{CaoBas_22}
X.~Cao and T.~Başar, ``Distributed constrained online convex optimization over multiple access fading channels,'' {\em IEEE Transactions on Signal Processing}, 2022.

\bibitem{GuoCaoHe_23}
H.~Guo, H.~Cao, J.~He, X.~Liu, and Y.~Shi, ``Pobo: Safe and optimal resource management for cloud microservices,'' {\em Performance Evaluation}, 2023.

\bibitem{CaoZhaPor_18}
X.~Cao, J.~Zhang, and H.~V. Poor, ``A virtual-queue-based algorithm for constrained online convex optimization with applications to data center resource allocation,'' {\em IEEE Journal of Selected Topics in Signal Processing}, 2018.

\bibitem{HsuXuLin_22}
W.-K. Hsu, J.~Xu, X.~Lin, and M.~R. Bell, ``Integrated online learning and adaptive control in queueing systems with uncertain payoffs,'' {\em Operations Research}, vol.~70, no.~2, pp.~1166--1181, 2022.

\bibitem{JohKamKan_21}
R.~Johari, V.~Kamble, and Y.~Kanoria, ``Matching while learning,'' {\em Operations Research}, vol.~69, no.~2, pp.~655--681, 2021.

\bibitem{McmHolScu_13}
H.~B. McMahan, G.~Holt, D.~Sculley, M.~Young, D.~Ebner, J.~Grady, L.~Nie, T.~Phillips, E.~Davydov, D.~Golovin, S.~Chikkerur, D.~Liu, M.~Wattenberg, A.~M. Hrafnkelsson, T.~Boulos, and J.~Kubica, ``Ad click prediction: A view from the trenches,'' KDD, 2013.

\bibitem{BalLuMir_20}
S.~Balseiro, H.~Lu, and V.~Mirrokni, ``Dual mirror descent for online allocation problems,'' in {\em Proceedings of the 37th International Conference on Machine Learning}, PMLR, 2020.

\bibitem{ZhaKeyEsm_20}
M.~Zhalechian, E.~Keyvanshokooh, C.~Shi, and M.~P. Van~Oyen, ``Personalized hospital admission control: A contextual learning approach,'' {\em Available at SSRN 3653433}, 2020.

\bibitem{TewSus_17}
A.~Tewari and S.~A. Murphy, {\em From Ads to Interventions: Contextual Bandits in Mobile Health}.
\newblock Springer International Publishing, 2017.

\bibitem{ChaKal_22}
S.~Chaudhary and D.~Kalathil, ``Safe online convex optimization with unknown linear safety constraints,'' {\em Proceedings of the AAAI Conference on Artificial Intelligence}, 2022.

\bibitem{ChaChaKal_23}
T.-J. Chang, S.~Chaudhary, D.~Kalathil, and S.~Shahrampour, ``Dynamic regret analysis of safe distributed online optimization for convex and non-convex problems,'' {\em Transactions on Machine Learning Research}, 2023.

\bibitem{AmaAliThr_19}
S.~Amani, M.~Alizadeh, and C.~Thrampoulidis, ``Linear stochastic bandits under safety constraints,'' in {\em Advances in Neural Information Processing Systems}, 2019.

\bibitem{AmaAliThr_20a}
A.~Moradipari, M.~Alizadeh, and C.~Thrampoulidis, ``Linear thompson sampling under unknown linear constraints,'' in {\em ICASSP 2020 - 2020 IEEE International Conference on Acoustics, Speech and Signal Processing (ICASSP)}, 2020.

\bibitem{AmaAliThr_20b}
S.~Amani, M.~Alizadeh, and C.~Thrampoulidis, ``Generalized linear bandits with safety constraints,'' in {\em ICASSP 2020 - 2020 IEEE International Conference on Acoustics, Speech and Signal Processing (ICASSP)}, 2020.

\bibitem{AmaAliThr_22}
A.~Moradipari, S.~Amani, M.~Alizadeh, and C.~Thrampoulidis, ``Safe linear thompson sampling with side information,'' {\em IEEE Transactions on Signal Processing}, 2021.

\bibitem{HutCheAli_24}
S.~Hutchinson, T.~Chen, and M.~Alizadeh, ``Optimistic safety for online convex optimization with unknown linear constraints,'' 2024.

\bibitem{MahJinYan_12}
M.~Mahdavi, R.~Jin, and T.~Yang, ``Trading regret for efficiency: online convex optimization with long term constraints,'' {\em The Journal of Machine Learning Research}, vol.~13, no.~1, pp.~2503--2528, 2012.

\bibitem{SunDeyKap_17}
W.~Sun, D.~Dey, and A.~Kapoor, ``Safety-aware algorithms for adversarial contextual bandit,'' in {\em International Conference on Machine Learning}, pp.~3280--3288, PMLR, 2017.

\bibitem{YuNeeWei_17}
H.~Yu, M.~Neely, and X.~Wei, ``Online convex optimization with stochastic constraints,'' {\em Advances in Neural Information Processing Systems}, vol.~30, 2017.

\bibitem{NeeYu_17}
M.~J. Neely and H.~Yu, ``Online convex optimization with time-varying constraints,'' {\em arXiv preprint arXiv:1702.04783}, 2017.

\bibitem{ChaSal_19}
N.~Liakopoulos, A.~Destounis, G.~Paschos, T.~Spyropoulos, and P.~Mertikopoulos, ``Cautious regret minimization: Online optimization with long-term budget constraints,'' in {\em Proceedings of the 36th International Conference on Machine Learning}, 2019.

\bibitem{CaoZhaPoo_21}
X.~Cao, J.~Zhang, and H.~V. Poor, ``Online stochastic optimization with time-varying distributions,'' {\em IEEE Transactions on Automatic Control}, vol.~66, no.~4, pp.~1840--1847, 2021.

\bibitem{SadRauFaz_20}
O.~Sadeghi, P.~Raut, and M.~Fazel, ``A single recipe for online submodular maximization with adversarial or stochastic constraints,'' in {\em Advances in Neural Information Processing Systems}, 2020.

\bibitem{YuaLam_18}
J.~Yuan and A.~Lamperski, ``Online convex optimization for cumulative constraints,'' {\em Advances in Neural Information Processing Systems}, vol.~31, 2018.

\bibitem{YiLiYan_21a}
X.~Yi, X.~Li, T.~Yang, L.~Xie, T.~Chai, and K.~Johansson, ``Regret and cumulative constraint violation analysis for online convex optimization with long term constraints,'' in {\em International Conference on Machine Learning}, pp.~11998--12008, PMLR, 2021.

\bibitem{GuoLiuWei_22}
H.~Guo, X.~Liu, H.~Wei, and L.~Ying, ``Online convex optimization with hard constraints: Towards the best of two worlds and beyond,'' in {\em Advances in Neural Information Processing Systems}, 2022.

\bibitem{LiuYanYin_23}
X.~Liu, Z.~Yang, and L.~Ying, ``Online nonstochastic control with adversarial and static constraints,'' in {\em Proceedings of the 40th International Conference on Machine Learning}, 2023.

\bibitem{AbhRah_24}
A.~Sinha and R.~Vaze, ``Optimal algorithms for online convex optimization with adversarial constraints,'' in {\em Advances in Neural Information Processing Systems}, 2024.

\bibitem{BadKleSli_13}
A.~Badanidiyuru, R.~Kleinberg, and A.~Slivkins, ``Bandits with knapsacks,'' in {\em Proceedings of the 2013 IEEE 54th Annual Symposium on Foundations of Computer Science}, FOCS '13, 2013.

\bibitem{BadKleSli_18}
A.~Badanidiyuru, R.~Kleinberg, and A.~Slivkins, ``Bandits with knapsacks,'' {\em Journal of the ACM (JACM)}, vol.~65, no.~3, pp.~1--55, 2018.

\bibitem{BadLanSli_14}
A.~Badanidiyuru, J.~Langford, and A.~Slivkins, ``Resourceful contextual bandits,'' in {\em Conference on Learning Theory}, pp.~1109--1134, PMLR, 2014.

\bibitem{AgrDev_16}
S.~Agrawal and N.~Devanur, ``Linear contextual bandits with knapsacks,'' {\em Advances in Neural Information Processing Systems}, vol.~29, 2016.

\bibitem{ImmSanSch_22}
N.~Immorlica, K.~Sankararaman, R.~Schapire, and A.~Slivkins, ``Adversarial bandits with knapsacks,'' {\em Journal of the ACM}, vol.~69, no.~6, pp.~1--47, 2022.

\bibitem{SivZuoBan_22}
V.~Sivakumar, S.~Zuo, and A.~Banerjee, ``Smoothed adversarial linear contextual bandits with knapsacks,'' in {\em Proceedings of the 39th International Conference on Machine Learning}, 2022.

\bibitem{BerCasCel_24b}
M.~Bernasconi, M.~Castiglioni, A.~Celli, and F.~Fusco, ``Bandits with replenishable knapsacks: the best of both worlds,'' in {\em The Twelfth International Conference on Learning Representations}, 2024.

\bibitem{CheCueLuo_20}
Y.~Chen, A.~Cuellar, H.~Luo, J.~Modi, H.~Nemlekar, and S.~Nikolaidis, ``Fair contextual multi-armed bandits: Theory and experiments,'' in {\em Proceedings of the 36th Conference on Uncertainty in Artificial Intelligence (UAI)}, 2020.

\bibitem{Sin_24}
A.~Sinha, ``Banditq: Fair bandits with guaranteed rewards,'' in {\em The 40th Conference on Uncertainty in Artificial Intelligence}, 2024.

\bibitem{SliSanFos_23}
A.~Slivkins, K.~A. Sankararaman, and D.~J. Foster, ``Contextual bandits with packing and covering constraints: A modular lagrangian approach via regression,'' in {\em Annual Conference Computational Learning Theory}, 2023.

\bibitem{SliSanFos_24}
A.~Slivkins, X.~Zhou, K.~A. Sankararaman, and D.~J. Foster, ``Contextual bandits with packing and covering constraints: A modular lagrangian approach via regression,'' {\em Journal of Machine Learning Research}, 2024.

\bibitem{GuoLiu_24}
H.~Guo and X.~Liu, ``Stochastic constrained contextual bandits via lyapunov optimization based estimation to decision framework,'' in {\em Proceedings of Thirty Seventh Conference on Learning Theory}, Proceedings of Machine Learning Research, 2024.

\bibitem{GuoLiu_25}
H.~Guo and X.~Liu, ``On stochastic contextual bandits with knapsacks in small budget regime,'' in {\em The Thirteenth International Conference on Learning Representations}, 2025.

\bibitem{PacGhaBar_21}
A.~Pacchiano, M.~Ghavamzadeh, P.~Bartlett, and H.~Jiang, ``Stochastic bandits with linear constraints,'' in {\em Proceedings of The 24th International Conference on Artificial Intelligence and Statistics}, 2021.

\bibitem{CheGanSal_22}
T.~Chen, A.~Gangrade, and V.~Saligrama, ``Strategies for safe multi-armed bandits with logarithmic regret and risk,'' in {\em Proceedings of the 39th International Conference on Machine Learning}, Proceedings of Machine Learning Research, 2022.

\bibitem{GanCheSal_24}
A.~Gangrade, T.~Chen, and V.~Saligrama, ``Safe linear bandits over unknown polytopes,'' in {\em Proceedings of Thirty Seventh Conference on Learning Theory}, 2024.

\bibitem{GanGopSal_24}
A.~Gangrade, A.~Gopalan, V.~Saligrama, and C.~Scott, ``Testing the feasibility of linear programs with bandit feedback,'' in {\em Proceedings of the 41st International Conference on Machine Learning}, Proceedings of Machine Learning Research, 2024.

\bibitem{AueCesFis_02}
P.~Auer, N.~Cesa-Bianchi, and P.~Fischer, ``Finite-time analysis of the multiarmed bandit problem,'' {\em Mach. Learn.}, 2002.

\bibitem{AbbPalSze_11}
Y.~Abbasi-yadkori, D.~P\'{a}l, and C.~Szepesv\'{a}ri, ``Improved algorithms for linear stochastic bandits,'' in {\em Advances in Neural Information Processing Systems 24}, 2011.

\bibitem{LatSze_20}
T.~Lattimore and C.~Szepesvári, {\em Bandit Algorithms}.
\newblock Cambridge University Press, 2020.

\bibitem{Nee_10}
M.~J. Neely, ``Stochastic network optimization with application to communication and queueing systems,'' {\em Synthesis Lectures on Communication Networks}, vol.~3, no.~1, pp.~1--211, 2010.

\bibitem{SriYin_14}
R.~Srikant and L.~Ying, {\em Communication Networks: {A}n Optimization, Control and Stochastic Networks Perspective}.
\newblock Cambridge University Press, 2014.

\bibitem{BoyVan_04}
S.~Boyd and L.~Vandenberghe, {\em Convex optimization}.
\newblock Cambridge university press, 2004.

\bibitem{SchRouBac_11}
M.~Schmidt, N.~Roux, and F.~Bach, ``Convergence rates of inexact proximal-gradient methods for convex optimization,'' in {\em Advances in Neural Information Processing Systems}, vol.~24, 2011.

\bibitem{Zin_03}
M.~Zinkevich, ``Online convex programming and generalized infinitesimal gradient ascent,'' in {\em Proceedings of the 20th International Conference on Machine Learning}, pp.~928--936, 2003.

\bibitem{HazKal_12}
E.~Hazan and S.~Kale, ``Projection-free online learning,'' in {\em Proceedings of the 29th International Coference on International Conference on Machine Learning}, 2012.

\bibitem{LiuLiShiYin_21}
X.~Liu, B.~Li, P.~Shi, and L.~Ying, ``An efficient pessimistic-optimistic algorithm for stochastic linear bandits with general constraints,'' in {\em Advances in Neural Information Processing Systems}, 2021.

\bibitem{GonMar_93}
G.~Chen and M.~Teboulle, ``Convergence analysis of a proximal-like minimization algorithm using bregman functions,'' {\em SIAM Journal on Optimization}, vol.~3, no.~3, pp.~538--543, 1993.

\bibitem{YuNee_20}
H.~Yu and M.~J. Neely, ``A low complexity algorithm with ${O}(\sqrt{T})$ regret and ${O}(1)$ constraint violations for online convex optimization with long term constraints,'' {\em The Journal of Machine Learning Research}, vol.~21, no.~1, pp.~1--24, 2020.

\bibitem{Lin_22}
L.~Xiao, ``On the convergence rates of policy gradient methods,'' {\em arXiv preprint arXiv:2201.07443}, 2022.

\bibitem{Har_13}
M.~Harchol-Balter, {\em Performance Modeling and Design of Computer Systems: Queueing Theory in Action}.
\newblock Cambridge University Press, 2013.

\bibitem{NewYorISO}
``New york {ISO} open access pricing data.,'' {\em http://www.nyiso.com/.}

\end{thebibliography}

\end{document}